\documentclass[12pt]{amsart}
\usepackage{fullpage}
\usepackage{times}
\usepackage{latexsym}
\usepackage{amssymb}
\usepackage[all]{xy}

\newtheorem{thm}{Theorem}[section]
\newtheorem{prop}[thm]{Proposition}
\newtheorem{cor}[thm]{Corollary}

\theoremstyle{remark}
\newtheorem{rem}[thm]{Remark}

\newtheorem{ex}[thm]{Example}

\usepackage{hyperref}

\newcommand{\CP}{\mathbb{CP}}

\newcommand{\Q}{\mathbb{Q}}
\newcommand{\R}{\mathbb{R}}
\newcommand{\Z}{\mathbb{Z}}

\title[Degrees of self-maps of products]
      {Degrees of self-maps of products}

\author{Christoforos Neofytidis}
\address{Department of Mathematical Sciences, SUNY Binghamton, Binghamton, NY 13902-6000, USA}
\email{chrisneo@math.binghamton.edu}
\date{\today; \copyright{ C. Neofytidis 2015}}

\begin{document}

\begin{abstract}
Every closed oriented manifold $M$ is associated with a set of integers $D(M)$, the set of self-mapping degrees of $M$. 
In this paper we investigate whether a product $M\times N$ admits a self-map of degree $d$, when neither $D(M)$ nor $D(N)$ contains $d$. We find sufficient conditions so that 
$D(M\times N)$ contains exactly the products of the elements of $D(M)$ with the elements of $D(N)$. As a consequence, we obtain manifolds $M\times N$ that do not admit self-maps of degree $-1$ (strongly chiral), that have finite sets of self-mapping degrees (inflexible) and that do not admit any self-map of degree $dp$ for a prime number $p$. Furthermore we obtain a characterization of odd-dimensional strongly chiral hyperbolic manifolds in terms of self-mapping degrees of their products. 
\end{abstract}

\maketitle

\section{Introduction}

Let $M$ be a closed oriented $n$-dimensional manifold and denote by $[M] \in H_n(M;\Z)$ its fundamental class. A continuous map $f \colon M \longrightarrow M$ is said to be of degree $d\in\Z$ if $H_n(f)([M])=d\cdot[M]$. The set of degrees of self-maps of $M$ is defined to be
\[
D(M):=\{d\in\Z \ | \ \exists \ f \colon M \longrightarrow M, \ \deg(f)=d\}.
\]

The investigation of $D(M)$ is a classical topic which has been used in several contexts to extract information about 
$M$, revealing simultaneously an interesting interplay between topology, global analysis and number theory. Obviously, every manifold satisfies $\{0,1\}\subseteq D(M)$.
In dimensions 1 and 2, all sets $D(M)$ are completely determined. In dimension 3, Wang et al. computed the unbounded sets $D(M)$, mostly following Thurston's geometrization picture; see~\cite{Wang,SWWZ,Waky} and the references given there. The set of self-mapping degrees of $M$ is bounded if and only if $M$ does not admit a self-map of absolute degree greater than one, that is $D(M)\subseteq\{-1,0,1\}$. A manifold with that property is termed {\em inflexible} in~\cite{CL}. So, the only remaining question in dimension $3$ is whether an inflexible manifold admits a self-map of degree $-1$. Prominent examples of inflexible manifolds (in any dimension) are the hyperbolic ones, because hyperbolic manifolds have positive simplicial volume~\cite{Gr1}. 
We refer to~\cite{Waky} for a concise description of the state of the art about $D(M)$ in dimension 3. In higher dimensions, Duan and Wang~\cite{DW} found a notable criterion, using the intersection form, for the (non-)existence of (self-)mapping degrees between highly connected manifolds of even dimension. 

For any two closed oriented manifolds $M$ and $N$, a trivial set-theoretic relation between $D(M)$, $D(N)$ and $D(M\times N)$ is given by
\begin{align}\label{setrelations}
D(M)\cup D(N) \subseteq D(M)\cdot D(N) \subseteq D(M\times N),
\end{align}
where $D(M)\cdot D(N):=\{\kappa\cdot\lambda \ |  \  \kappa\in D(M), \ \lambda\in D(N)\}$.
However, the converse inclusions do not generally hold; we illustrate this with two examples:

\begin{ex}\label{ex:-l2}
Let $N=\Sigma$ be a closed oriented surface of genus at least $2$. Since $\Sigma$ is hyperbolic, we conclude that $D(\Sigma)\subseteq\{-1,0,1\}$. In addition, $\Sigma$ admits a self-map of degree $-1$, being a connected sum of tori. Thus $D(\Sigma)=\{-1,0,1\}$. 
Now, let $M$ be a torus bundle over $S^1$ with monodromy 
$\left( \begin{array}{cc}
1 & 0 \\
a & 1 \end{array} \right), \ a\neq0$.
The set of self-mapping degrees of $M$ is $D(M)=\{k^2 \ | \ k\in\Z\}$; cf.~\cite{SWWZ}.
Thus, $-k^2\in D(M)\cdot D(\Sigma)$ for every $k\in\Z$. However, $-k^2\notin D(M)\cup D(\Sigma)$ for $k\neq0,1$, and so the converse of the first inclusion in (\ref{setrelations}) fails. 
\end{ex}

\begin{ex}\label{ex:-1}
Let $M$ be a closed oriented manifold of odd dimension. The map $(x,y)\mapsto(y,x)$ is an orientation-reversing involution of $M\times M$, which means that $-1\in D(M\times M)$. If we pick $M$ such that $-1\notin D(M)$, then $-1\notin D(M)\cdot D(M)$, and so the converse of the second inclusion in (\ref{setrelations}) fails. 
We note that there exists a variety of examples of manifolds that do not admit self-maps of degree $-1$. We refer to a result of Belolipetsky and Lubotzky~\cite{BL} (cf. Theorem \ref{t:BL} below) to deduce the existence of hyperbolic manifolds that do not admit self-maps of degree $-1$ in every dimension $\geq3$. Of course, there exist many other simpler examples of manifolds with that property; 
see the discussion in Section \ref{s:examples}. 
\end{ex}

Example \ref{ex:-1} originates from the following general question:

\medskip
{\em Suppose $M$ and $N$ are two closed oriented manifolds which have a (topological, diffeomorphism, etc.) property $P$. Does $M\times N$ have the property $P$ as well?}
\medskip

For instance, if the property $P$ is ``$-1\notin D(M)$", then Example \ref{ex:-1} says that although a manifold $M$ can have $P$, the direct product $M\times M$ does not necessarily have $P$ (see also Proposition \ref{p:hyper}).

\begin{rem}
The latter source of examples, concerning maps of degree $-1$, has close connections to corresponding concepts of mathematical biology and chemistry, most notably the notions of ``chiral knots" and ``chiral molecules". Due to this relation, a manifold $M$ satisfying $-1 \notin D(M)$ is called {\em strongly chiral}; we refer to~\cite{Mue} and the related references there for more details. 
\end{rem}

The following list of problems from~\cite[Section 1]{Sasao} and~\cite[Section 9: Appendix II]{CL} (see also~\cite[Section 5.35]{Gromov}) is the main motivation for this paper:

\begin{itemize}
\item[(i)] When does it hold $D(M)\cup D(N) = D(M)\cdot D(N) = D(M\times N)$?
\item[(ii)] Characterize $M$ so that $-1 \notin D(M)$, that is $M$ is strongly chiral. 
\item[(iii)] For a prime number $p$, characterize $M$ so that $p\Z\cap D(M)=\{0\}$.
\item[(iv)] When is the product of inflexible manifolds inflexible itself?
\end{itemize}

Our goal is to investigate how the existence of a mapping degree in $D(M\times N)$ reflects on the individual sets $D(M)$ and $D(N)$.   
We have seen in Example \ref{ex:-1} that neither $D(M)\cup D(N)$ nor $D(M)\cdot D(N)$ are in general optimal approximations for $D(M\times N)$, and Example \ref{ex:-l2} suggests that $D(M)\cdot D(N)$ is a better approximation. Our main result is that $D(M)\cdot D(N)$ is indeed equal to $D(M\times N)$ in certain cases:

\begin{thm}\label{t:main}
Let $M$ and $N$ be two closed oriented manifolds of dimensions $m$ and $n$ respectively. Suppose that there is no map of non-zero degree from a non-trivial direct product to $M$ and that $M$ cannot be realized by a cohomology class in $H^m(N;\Q)$\footnote{That is, the cohomological fundamental class of $M$ cannot be pulled back to a non-trivial class in $H^m(N;\Q)$ under a continuous map $N\longrightarrow M$.}. 
Then $D(M\times N)=D(M)\cdot D(N)$.
\end{thm}

The assumptions of Theorem \ref{t:main} occur naturally quite often. On the one hand, there are plenty of examples of manifolds that do not admit maps of non-zero degree from direct products~\cite{KL,KL2,KN,NeoThesis,NeoIIPP}. On the other hand, the assumption that $M$ cannot be realized by a cohomology class in $H^m(N;\Q)$ is fulfilled in several instances, e.g. when there is no map of non-zero degree from $N$ to $M$ (if $N$ has the same dimension as $M$), or, simply, when the dimension of $N$ is smaller than the dimension of $M$ (and so $H^m(N)$ is trivial); see Section \ref{s:examples} for examples.

For instance, the manifolds $M$ and $\Sigma$ in Example \ref{ex:-l2} fulfill the assumptions of Theorem \ref{t:main}, because $M$ does not admit maps of non-zero degree from direct products~\cite{KN} and $H^3(\Sigma)=0$. Thus, Theorem \ref{t:main} implies that $D(M\times\Sigma)=D(M)\cdot D(\Sigma)$. Since $D(M)\cup D(\Sigma) \subsetneq D(M)\cdot D(\Sigma)$, we conclude that $D(M)\cdot D(N)$ is indeed a considerably better approximation of $D(M\times N)$ than $D(M)\cup D(N)$. Moreover, Example \ref{ex:-l2} shows that the equality $D(M)\cdot D(N) = D(M\times N)$ in Theorem \ref{t:main} is the best we can obtain regarding Problem (i). 
Nevertheless, applying Theorem \ref{t:main} to certain mapping degrees in $D(M\times N)$, we can show that those degrees belong not only to $D(M)\cdot D(N)$ but also to $D(M)\cup D(N)$. More precisely, we obtain the following characterizations with respect to Problems (ii) -- (iv):

\begin{cor}\label{c:main}
Let $M$ and $N$ be as in Theorem \ref{t:main}. Then the following hold:
\begin{itemize}
\item[(a)] $-1\notin D(M)\cup D(N)$ if and only if $-1\notin D(M\times N)$. That is, $M\times N$ is strongly chiral if and only if both $M$ and $N$ are strongly chiral.
\item[(b)] For any prime number $p$, \[p\Z\cap(D(M)\cup D(N))=\{0\} \ \text{ if and only if }\ p\Z\cap D(M\times N)=\{0\}.\]
\item[(c)] $D(M\times N)\subseteq\{-1,0,1\}$ if and only if $D(M)\cup D(N)\subseteq\{-1,0,1\}$. That is, $M\times N$ is inflexible if and only if both $M$ and $N$ are inflexible. 
\end{itemize}
\end{cor}

As we have seen in Example \ref{ex:-1}, part (a) of the above corollary does not hold when the requirement that $M$ cannot be realized by a class in $H^m(N;\Q)$ is violated. In fact, Mostow's rigidity implies that the phenomenon of Example \ref{ex:-1} characterizes up to isometry strongly chiral hyperbolic manifolds of the same odd dimension.

\begin{prop}\label{p:hyper}
Let $M$ and $N$ be two closed oriented hyperbolic manifolds. Then $D(M\times N)=D(M)\cup D(N)$, unless $M$ and $N$ are isometric and have odd dimension, in which case we have $D(M\times N)=D(M)\cup D(N)$ if and only if $-1\in D(M)\cup D(N)$.
\end{prop}

In particular, the product of two strongly chiral hyperbolic manifolds of the same odd dimension is strongly chiral if and only if those manifolds are not isometric. Examples of strongly chiral products of hyperbolic manifolds exist when the dimensions of the factors are at least three; see Corollary \ref{c:hyper}.

\subsection*{Acknowledgements} 
We would like to thank Shicheng Wang for a stimulating question about the relation between the sets $D(M)\cdot D(N)$ and $D(M\times N)$. We also thank the anonymous referees for useful comments on the structure of this paper.

\section{Proofs}

\subsection{Realization of (co-)homology classes by closed manifolds} 

One of the basic ingredients of the proofs is the following theorem of Thom~\cite{Thom}, which answers in the affirmative (in rational homology) Steenrod's classical problem~\cite[Problem 25]{Steenrodproblem} of the realization of (integral) homology classes by closed manifolds:

\begin{thm}[Thom's Realization Theorem~\cite{Thom}]\label{t:Thom}
Let $X$ be a topological space. For every $w\in H_n(X; \Z)$ there is an integer $d>0$ and a closed oriented smooth $n$-dimensional manifold $M$ together with a continuous map $f\colon M\longrightarrow X$ so that $H_n(f)([M])=d\cdot w$. 
In particular, every rational homology class in degree $n$ is realizable by a closed oriented smooth $n$-dimensional manifold.
\end{thm}

In the proof of Theorem \ref{t:main}, we will use the dual version of Thom's theorem in cohomology. Namely, if $\alpha\in H^n(X;\Z)$, then Theorem \ref{t:Thom} states that there exists an integer $d>0$ and a closed oriented smooth $n$-dimensional manifold $M$ together with a continuous map $f\colon M\longrightarrow X$ so that $H^n(f)(\alpha)=d\cdot\omega_M$, where $\omega_M\in H^n(M;\Z)$ is the cohomological fundamental class of $M$.

\subsection{Proof of Theorem \ref{t:main}} 
It suffices to show that
\[
D(M\times N)\subseteq D(M)\cdot D(N).
\]
Suppose $f\colon M\times N\longrightarrow M\times N$ is a map of degree $d$. The K\"unneth formula for the cohomology group (with rational coefficients) of $M\times N$ in degree $l\in\{0,...,m+n\}$ is
\begin{align}\label{Kuenneth}
H^l(M\times N)\cong H^l(M) \oplus(H^{l-1}(M)\otimes H^1(N))\oplus\cdots\oplus H^{l}(N).
\end{align}
Let $p_M\colon M\times N\longrightarrow M$ be the projection onto $M$ and consider the composite map
\[
M\times N\stackrel{f}\longrightarrow M\times N\stackrel{p_M}\longrightarrow M.
\]
Since $M$ does not admit maps of non-zero degree from direct products, Theorem \ref{t:Thom} implies that $\omega_M$ maps trivially under $H^m(p_M\circ f)$ to all products $H^i(M)\otimes H^j(N)$, where $0< i,j <m$ and $i+j=m$ (see also~\cite{KLN}). Moreover, since $\omega_M$ cannot be realized by any class $\beta_m\in H^m(N;\Q)$, we conclude that
\begin{align}\label{omegaM}
H^m(f)(\omega_M\times 1)=\kappa\cdot(\omega_M\times 1)\in H^m(M\times N),
\end{align}
for some $\kappa\in\Z$.

Similarly, by the K\"unneth formula (\ref{Kuenneth}) in degree $n$ we obtain
\begin{align}\label{omegaN}
H^n(f)(1\times\omega_N)=\sum_{i=0}^{n}\lambda_i\cdot(\alpha_{n-i}\times\beta_i)\ \in H^n(M\times N),
\end{align}
where $\lambda_i\in\Z, \ \alpha_{n-i}\in H^{n-i}(M)$ and $\beta_i\in H^i(N)$, for $i=0,1,...,n$; in particular, $\beta_n=\omega_N$. 

Combining (\ref{omegaM}) and (\ref{omegaN}), and by the naturality of the cup product, we obtain
\[
d\cdot\omega_{M\times N}=H^{m+n}(f)(\omega_{M\times N})=\kappa\lambda_n\cdot\omega_{M\times N}.
\]
Let $\iota_M\colon M\hookrightarrow M\times N$ be the inclusion of $M$ into $M\times N$. Then $\kappa$ is realized as the degree of the composite map
\[
M\stackrel{\iota_M}\hookrightarrow M\times N\stackrel{f}\longrightarrow M\times N\stackrel{p_M}\longrightarrow M.
\]
Similarly, $\lambda_n$ is realized as the degree of the map $p_N\circ f\circ\iota_N\colon N\longrightarrow N$, where again $\iota_N$ and $p_N$ denote inclusion and projection respectively. This means that $d\in D(M)\cdot D(N)$ as desired.

\subsection{Proof of Corollary \ref{c:main}}
\begin{itemize}
\item[(a)] It suffices to show that $-1\in D(M\times N)$ implies $-1\in D(M)\cup D(N)$. If $-1\in D(M\times N)$, then $-1\in D(M)\cdot D(N)$ by Theorem \ref{t:main}, and so clearly $-1\in D(M)\cup D(N)$.
\item[(b)] Let $p$ be a prime number. We will show that, if there exists a non-zero integer $\ell$ so that $p\cdot\ell\in D(M\times N)$, then $p\Z\cap(D(M)\cup D(N))\neq\{0\}$. For $\deg(f)=p\cdot\ell$ in Theorem \ref{t:main}, we have that $p\cdot\ell\in D(M)\cdot D(N)$. Since $p$ is prime, we conclude that $p$ divides at least one mapping degree in $D(M)$ or $D(N)$, which means that $p\Z\cap(D(M)\cup D(N))\neq\{0\}$ as required.
\item[(c)] We only need to show that $M\times N$ is inflexible, when both $M$ and $N$ are inflexible. If $f\colon M\times N\longrightarrow M\times N$ is a map of degree $d$, then $d\in D(M)\cdot D(N)$ by Theorem \ref{t:main}. Since $D(M), D(N)\subseteq\{-1,0,1\}$, we conclude that $d\in\{-1,0,1\}$.
\end{itemize}

\subsection{Proof of Proposition \ref{p:hyper}}
Since $D(M)\cup D(N)\subseteq D(M\times N)\subseteq\{-1,0,1\}$ (recall that $M\times N$ has also positive simplicial volume~\cite{Gr1}), it suffices to examine what happens when $-1\in D(M\times N)$. Let $\dim (M)=m$ and $\dim (N)=n$ and suppose $f\colon M\times N\longrightarrow M\times N$ is a map of degree $-1$. By a result of Kotschick and L\"oh~\cite{KL}, hyperbolic manifolds do not admit maps of non-zero degree from direct products. Thus, Theorem \ref{t:Thom} implies that
\begin{align}\label{eq:hyper1}
H^m(f)(\omega_M\times 1)=\kappa\cdot(\omega_M\times 1)+\mu\cdot (1\times\beta_m), \ \beta_m\in H^m(N)
\end{align}
and
\begin{align}\label{eq:hyper2}
H^n(f)(1\times\omega_N)=\nu\cdot(\alpha_n\times 1)+\lambda\cdot (1\times\omega_N), \ \alpha_n\in H^n(M)
\end{align}
for some integers $\kappa,\mu,\nu,\lambda$; compare (\ref{omegaM}) and (\ref{omegaN}) respectively. Also, note that $\kappa\in D(M)$ and $\lambda\in D(N)$, which means that $\kappa,\lambda\in\{-1,0,1\}$.

If $m\neq n$, then $\beta_m=0$ or $\alpha_n=0$, and so (\ref{eq:hyper1}) and (\ref{eq:hyper2}) imply that $-1=\deg(f)=\kappa\lambda$. This means that $-1\in D(M)\cup D(N)$.

Suppose now $m=n$. In this case, $\beta_m=\omega_N$ in (\ref{eq:hyper1}) and $\alpha_n=\omega_M$ in (\ref{eq:hyper2}). Thus
\begin{align}\label{degree}
-1=\deg(f)=\kappa\lambda+(-1)^{m}\mu\nu
\end{align}
Note that $\mu$ and $\nu$ can be realized as degrees for maps $N\longrightarrow M$ and $M\longrightarrow N$ respectively (through the maps $p_M\circ f\circ \iota_N$ and $p_N\circ f\circ \iota_M$ respectively, where $\iota$ denotes inclusion and $p$ projection). If $\mu\nu=0$, then (\ref{degree}) becomes $\kappa\lambda=-1$, and so $-1\in D(M)\cup D(N)$. Suppose now $\mu\nu\neq0$. Then $\mu,\nu\in\{\pm1\}$, because $M$ and $N$ have non-zero simplicial volume~\cite{Gr1}, and so (\ref{degree}) becomes 
\[
\mu\nu=
\begin{cases}
-1 &\text{if $m$ is even},\\
1 & \text{if $m$ is odd}.
\end{cases}
\]
If $m$ is even, then $-1\in D(M)\cup D(N)$. However, when $m$ is odd, the equation $\mu\nu=1$ does not imply that $-1\in D(M)\cup D(N)$. Also, note that $M$ and $N$ are isometric whenever $\mu\nu\neq0$: There are maps $g_1\colon M\longrightarrow N$ and $g_2\colon N\longrightarrow M$ of absolute degree one and the composite map $g_2\circ g_1\colon M\longrightarrow M$ is again of absolute degree one, which implies that the induced homomorphism $(g_2\circ g_1)_*\colon \pi_1(M)\longrightarrow\pi_1(M)$ is surjective. 
By a classical theorem of Mal'cev on linear groups (or by a result of Sela~\cite{Sela}), $\pi_1(M)$ is Hopfian, which means that $(g_2\circ g_1)_*$ is an isomorphism. Similarly, the map $g_1\circ g_2$ induces an automorphism of $\pi_1(N)$. We conclude that each of the $g_i$ induces an isomorphism between $\pi_1(M)$ and $\pi_1(N)$, and so it is a homotopy equivalence. Finally, Mostow's rigidity theorem implies that $M$ and $N$ are isometric. 

Thus, the only case where the equality $D(M\times N)=D(M)\cup D(N)$ might not occur is when $m$ is odd and $M$ and $N$ are isometric. Let $M$ and $N$ be isometric of odd dimension and $h\colon M\longrightarrow N$ be an isometry. Then an orientation-reversing self-isometry of $M\times N$ is given by the map
\begin{align*}
M\times N &  \longrightarrow   
\ M\times N\\
 (x,y) & \ \mapsto  \  
 (h^{-1}(y),h(x)).
\end{align*}
(Note that the hypothesis that $M$ and $N$ are hyperbolic is not necessary here.) Thus the equality $D(M\times N)= D(M)\cup D(N)$ holds if and only if $-1\in D(M)\cup D(N)$. 

The proof is now complete.

\section{Examples}\label{s:examples}

\subsection{Products, inflexibility and chirality}

The non-existence of maps of non-zero degree from direct products (of surfaces) to certain aspherical manifolds was raised by Gromov~\cite{Gromov} in his theory of bounded cohomology and topological rigidity. Obstructions to the existence of such maps were developed recently~\cite{KL,KL2,NeoIIPP, NeoThesis}. Prominent examples of manifolds that do not admit maps of non-zero degree from direct products are low-dimensional aspherical manifolds that possess a non-product Thurston geometry and manifolds with non-positive sectional curvature that are not virtual products themselves.
In particular, as we already mentioned above, no hyperbolic manifold admits a map of non-zero degree from a direct product~\cite{KL}. 

Moreover, we have seen that every hyperbolic manifold is inflexible, and so the remaining question is to determine which hyperbolic manifolds are strongly chiral. In dimension 2, every closed hyperbolic surface admits an orientation reversing self-diffeomorphism (see also Example \ref{ex:-l2}), however in higher dimensions there is not a complete answer. Nevertheless, the following result of Belolipetsky and Lubotzky~\cite{BL} implies the existence of hyperbolic manifolds that do not admit self-maps of degree $-1$ in every dimension $\geq 3$:

\begin{thm}[\cite{BL}]\label{t:BL}
For every $n\geq 2$ and every finite group $\Gamma$ there exist closed oriented $n$-dimensional hyperbolic manifolds $M$ with $\mathrm{Isom}(M)\cong\Gamma$.
\end{thm}

If $\Gamma$ is of odd order and $n\geq 3$, then the hyperbolic manifolds of the above theorem are strongly chiral by Mostow's rigidity theorem, as observed by Weinberger; see~\cite[Section 3]{Mue}. Hence, we have the following consequence of Theorem \ref{t:BL} and Proposition \ref{p:hyper}, providing examples of strongly chiral products of hyperbolic manifolds:

\begin{cor}\label{c:hyper}
For every $m,n\geq3$ there exist closed oriented strongly chiral hyperbolic manifolds $M$ and $N$ of dimensions $m$ and $n$ respectively, and their product $M\times N$ is always strongly chiral, unless $M$ is isometric to $N$ and $m$ is odd.
\end{cor}

Several other obstructions to the existence of self-maps of degree $-1$ were developed in the past, using, for instance, the intersection form in dimensions $4n$ and the linking form in dimensions $4n-1$. Among the most standard examples of strongly chiral manifolds are the complex projective spaces $\CP^{2n}$. M\"ullner~\cite{Mue} showed that in each dimension $\geq7$ there exist simply connected manifolds that do not admit self-maps of degree $-1$. 
Using $\CP^{2n}$ and certain $S^{2n-1}$-bundles over $S^{2n}$ as the main building factors, the proof given in~\cite{Mue} is based on the following:

\begin{prop}[\normalfont{\cite[Section 3]{Mue}}]\label{p:Mue}
Let $M$ be a rational homology sphere of dimension $m$ and $N$ be a closed oriented manifold of dimension $n$ such that either
\begin{itemize}
\item[(1)] $N$ is not a rational homology sphere, if $m=n$, or
\item[(2)] $H^m(N;\Q)=0$, if $m\neq n$.
\end{itemize}
Then $M\times N$ is strongly chiral if and only if both $M$ and $N$ are strongly chiral.
\end{prop}

The proof of above proposition relies on two facts: First, since $M$ is a rational homology sphere, the K\"unneth theorem gives $H^m(M\times N;\Q)\cong H^m(M;\Q)\oplus H^m(N;\Q)$.
This means that for any self-map $f\colon M\times N\longrightarrow M\times N$ there exist integers $\kappa,\mu$ such that $H^m(f)(\omega_M\times 1)=\kappa\cdot(\omega_M\times 1)+\mu\cdot(1\times\beta_m)$, where $\beta_m\in H^m(N;\Q)$. Second, either there is no map of non-zero degree from $M$ to $N$, when the assumption (1) holds\footnote{If $f\colon M\longrightarrow N$ is a non-zero degree map, then the induced homomorphisms 
$H^*(f)\colon H^*(N;\Q)\longrightarrow H^*(M;\Q)$ are injective.}, 
or $H^m(M\times N;\Q)\cong H^m(M;\Q)$ -- and so $M$ cannot be realized by any class $\beta_m\in H^m(N;\Q)$ -- when the assumption (2) is satisfied. Thus, Theorem \ref{t:main} can be viewed as a (co-)homological extension of the idea of Proposition \ref{p:Mue}, since it does not require anymore the vanishing of the product groups $H^j(M;\Q)\otimes H^{m-j}(N;\Q)$, $j\neq0,m$, or of $H^{m}(N;\Q)$.

\subsection{Applications}

One of our basic building factors will be hyperbolic manifolds, although we could more generally consider irreducible locally symmetric spaces of non-compact type~\cite{LS,KL}. 
Corollary \ref{c:hyper} gives already examples of strongly chiral products of hyperbolic manifolds $M$ and $N$ in all dimensions $\dim(M),\dim(N)\geq 3$.

Following Corollary \ref{c:main} (a), further examples (together with their generalizations) of strongly chiral products can be obtained as follows: 

\begin{ex}
As mentioned above, the complex projective plane $\CP^2$ is strongly chiral. 
Moreover, since $\CP^2$ is simply connected, any map from $\CP^2$ to a closed $4$-manifold $M$ with infinite fundamental group has degree zero. Let for example $M$ be a strongly chiral hyperbolic $4$-manifold.
Then Corollary \ref{c:main} (a) implies that $-1\notin D(M\times\CP^2)$, i.e. $M\times\CP^2$ is strongly chiral.

Observe that the conclusion that any map from $\CP^2$ to $M$ has degree zero can be deduced as well by the fact that $\CP^2$ admits a dominant map from the product $S^2\times S^2$ (a branched $2$-fold covering given as the quotient map of the involution $(x,y)\mapsto(y,x)$ of $S^2\times S^2$). This means that $\CP^2$ may be replaced by any strongly chiral, closed oriented $4$-manifold $N$ that admits a map of non-zero degree from a direct product. In the same spirit, $M$ can be any strongly chiral closed oriented $4$-manifold that does not admit dominant maps from products.
\end{ex}

Using the computations of~\cite{SWWZ}, we can obtain products whose sets of self-mapping degrees do not contain non-trivial multiples of $p=2$:

\begin{ex}
According to~\cite[Theorem 1.3]{SWWZ}, there exist tori semi-bundles $N_1$, $N_2$ and $N_3$ possessing the geometries $\R^3$, $Nil^3$ and $Sol^3$ respectively, so that 
\[
D(N_1)=\{0, 2k+1 \ | \ k\in\Z\} \ 
\text{and} \ D(N_2)=D(N_3)=\{0, (2k+1)^2  \ | \ k\in\Z\}.
\]
In particular, $2\Z\cap D(N_i)=\{0\}$ for all $i\in\{1,2,3\}$. 
The manifolds $N_2$ and $N_3$ do not admit maps of non-zero degree from direct products~\cite{Wang1,KN} (however, $N_1$ is finitely covered by $T^3$ by a classical result of Bieberbach).  Moreover, they do not exist maps of non-zero degree between $N_i$ and $N_j$ for $i\neq j$; see~\cite{Wang1,KN}. Thus, Corollary \ref{c:main} (b) implies that 
\[2\Z\cap D(N_i\times N_j)=\{0\} \ \text{for all}\  i\neq j.\]
\end{ex}
\medskip

Furthermore, we can combine Theorem \ref{t:main} with the results of~\cite{SWWZ} to compute the sets of self-mapping degrees for several classes of products $M\times N$, where $N$ is any closed oriented $3$-manifold that possesses a Thurston geometry, 
and $M$ is a suitable manifold that does not admit maps of non-zero degree from products. For example:

\begin{ex}
Let $M$ be a closed oriented hyperbolic manifold of dimension $m\geq 4$ and $N$ be a closed oriented $3$-manifold that possesses a Thurston geometry. The set $D(N)$ is either $\{0,1\}$, $\{-1,0,1\}$, or infinite and explicitly computed in~\cite{SWWZ}. Since $M$ does not admit dominant maps from products~\cite{KL} and $H^m(N)=0$ (because $m\geq4$), Theorem \ref{t:main} implies that $D(M\times N)=D(M)\cdot D(N)$. Thus
\[
D(M\times N)=
\begin{cases}
D(N) &\text{if $M$ is strongly chiral},\\
D(N)\cup(-D(N)) & \text{otherwise}.
\end{cases}
\]
\end{ex}

\medskip

Obstructions to the existence of self-maps of absolute degree greater than one can be derived by the positivity of numerical invariants $I\in[0,\infty]$ that are monotonous under continuous maps. That is, if $f\colon M\longrightarrow M$ is a map of non-zero degree, then $I(M)\geq |\deg(f)|\cdot I(M)$, which implies that $|\deg(f)|\leq1$ whenever $I(M)>0$. In this paper, our basic examples of inflexible manifolds (and one of the building factors for constructing products) were the hyperbolic ones, because hyperbolic manifolds have positive simplicial volume.
A product of two hyperbolic manifolds $M\times N$ is again inflexible, because the simplicial volume 
satisfies $\|M\times N\|\geq\|M\|\|N\|$; cf.~\cite{Gr1}. However, simply connected manifolds have zero simplicial volume and the same holds for all products containing a simply connected factor~\cite{Gr1}. In fact, it is an open question whether there is a finite semi-norm that does not vanish on a simply connected manifold~\cite{Gromov,CL}. Nevertheless, simply connected inflexible manifolds do exist, at least in high dimensions~\cite{CL,Am,CV}. Using those examples we can obtain inflexible products that contain at least one simply connected factor (and thus have vanishing simplicial volume):

\begin{ex}
Let $M$ be a closed oriented manifold with positive simplicial volume that does not admit maps of non-zero degree from direct products, 
and $N$ be an inflexible, closed oriented simply connected manifold of dimension $\leq\dim(M)$. Since $\pi_1(M)$ is infinite (because $\|M\|>0$), there is no map of non-zero degree from $N$ to $M$, and so Corollary \ref{c:main} (c) implies that $M\times N$ is inflexible.
\end{ex}

We remark that inflexible products of simply connected manifolds were given in~\cite[Section 9]{CL}.

\bibliographystyle{amsplain}

\end{document}